# DYNAMICS OF THE TIME TO THE MOST RECENT COMMON ANCESTOR IN A LARGE BRANCHING POPULATION

BY STEVEN N. EVANS[1] AND PETER L. RALPH[2]

*University of California at Berkeley and University of California at Davis*

If we follow an asexually reproducing population through time, then the amount of time that has passed since the most recent common ancestor (MRCA) of all current individuals lived will change as time progresses. The resulting "MRCA age" process has been studied previously when the population has a constant large size and evolves via the diffusion limit of standard Wright–Fisher dynamics. For any population model, the sample paths of the MRCA age process are made up of periods of linear upward drift with slope $+1$ punctuated by downward jumps. We build other Markov processes that have such paths from Poisson point processes on $\mathbb{R}_{++} \times \mathbb{R}_{++}$ with intensity measures of the form $\lambda \otimes \mu$ where $\lambda$ is Lebesgue measure, and $\mu$ (the "family lifetime measure") is an arbitrary, absolutely continuous measure satisfying $\mu((0, \infty)) = \infty$ and $\mu((x, \infty)) < \infty$ for all $x > 0$. Special cases of this construction describe the time evolution of the MRCA age in $(1 + \beta)$-stable continuous state branching processes conditioned on nonextinction—a particular case of which, $\beta = 1$, is Feller's continuous state branching process conditioned on nonextinction. As well as the continuous time process, we also consider the discrete time Markov chain that records the value of the continuous process just before and after its successive jumps. We find transition probabilities for both the continuous and discrete time processes, determine when these processes are transient and recurrent and compute stationary distributions when they exist. Moreover, we introduce a new family of Markov processes that stands in a relation with respect to the general $(1+\beta)$-stable continuous state branching process and its conditioned version that is similar to the one between

Received December 2008; revised April 2009.

[1]Supported in part by NSF Grants DMS-04-05778 and DMS-09-07630.

[2]Supported in part by a VIGRE grant awarded to the Department of Statistics, University of California at Berkeley.

the family of Bessel-squared diffusions and the unconditioned and conditioned Feller continuous state branching process.

**1. Introduction.** Any asexually reproducing population has a unique most recent common ancestor, from whom the entire population is descended. In sexually reproducing species, the same is true for each nonrecombining piece of DNA. For instance, our "mitochondrial Eve," from whom all modern-day humans inherited their mitochondrial DNA, is estimated to have lived around 180,000 years ago [19] while our "Y-chromosomal Adam" is estimated to have lived around 50,000 years ago [33]. There have also been efforts to estimate the time since the MRCA lived (which we will also call the "age of the MRCA") in populations of other organisms, particularly pathogens [32, 34]. These studies, using sophisticated models incorporating of demographic history, are focused on estimating the age of the MRCA at a single point in time (the present).

As time progresses into the future, eventually the mitochondrial lineages of all but one of the daughters of the current mitochondrial Eve will die out, at which point the new mitochondrial Eve will have lived somewhat later in time. The age of the MRCA is thus a dynamically evolving *process* that exhibits periods of upward linear growth separated by downward jumps.

Recently, [26] and [29] independently investigated the MRCA age process for the diffusion limit of the classical Wright–Fisher model. The Wright–Fisher model is perhaps the most commonly used model in population dynamics: each individual in a fixed size population independently gives birth to an identically distributed random number of individuals (with finite variance), and after the new offspring are produced, some are chosen at random to survive so that the total population size remains constant. The diffusion limit arises by letting the population size go to infinity and taking the time between generations to be proportional to the reciprocal of the population size.

In this paper, we investigate the MRCA age process for a parametric family of population models in a setting in which the population size varies with time, and, by suitable choice of parameters, allows control over the extent to which rare individuals can have large numbers of offspring that survive to maturity. The model for the dynamics of the population size is based on the *critical $(1+\beta)$-stable continuous state branching process* for $0 < \beta \leq 1$. These processes arise as scaling limits of Galton–Watson branching processes as follows.

Write $Z_t^{(n)}$ for the number of individuals alive in a critical continuous time Galton–Watson branching process with branching rate $\lambda$ and offspring distribution $\gamma$. The distribution $\gamma$ has mean 1 (and thus, the process is "critical"). Suppose that if $W$ is a random variable with distribution $\gamma$, then the random walk with steps distributed as the random variable $(W - 1)$ falls into the domain of attraction of a stable process of index $1 + \beta \in (1, 2]$. The case $\beta = 1$ corresponds to $\gamma$ having finite variance and the random walk converging to Brownian motion after rescaling. Set $X_t^{(n)} = n^{-1/\beta} Z_t^{(n)}$ and suppose that $X_0^{(n)} \to x$ as $n \to \infty$. Then, up to a time-rescaling depending on $\lambda$ and the scaling of the stable process above, the processes $X^{(n)}$ converge to a Markov process $X$ that is a critical $(1 + \beta)$-stable continuous state branching process, and whose distribution is determined by the Laplace transform

$$(1.1) \qquad \mathbb{E}[e^{-\theta X_t} | X_0 = x] = \exp\left(-\frac{\theta x}{(1 + \theta^\beta t)^{1/\beta}}\right).$$

If $\beta = 1$, this is Feller's critical continuous state branching process [17, 22]. (Note that time here is scaled by a factor of 2 relative to some other authors so that the generator of our "Feller continuous state branching process" is $x \frac{\partial^2}{\partial x^2}$.)

Let $\tau = \inf\{t > 0 : X_t = 0\}$ denote the *extinction time* of $X$ (it is not hard to show that $X_t = 0$ for all $t \geq \tau$). Taking $\theta \to \infty$ in (1.1) gives

$$\mathbb{P}\{\tau > t | X_0 = x\} = 1 - \exp\left(-\frac{x}{t^{1/\beta}}\right),$$

so $X$ dies out almost surely. However, it is possible to condition $X$ to live forever in the following sense:

$$\lim_{T \to \infty} \mathbb{E}[f(X_t) | X_0 = x, \tau > T] = \frac{1}{x} \mathbb{E}[f(X_t) X_t | X_0 = x].$$

Thus if $P_t(x', dx'')$ are the transition probabilities of $X$, then there is a Markov process $Y$ with transition probabilities

$$Q_t(y', dy'') = \frac{1}{y'} P_t(y', y'') y''.$$

The process $Y$ is the critical $(1 + \beta)$-stable continuous state branching process $X$ *conditioned on nonextinction*. The distribution of $Y_t$ is determined by its Laplace transform

$$(1.2) \quad \mathbb{E}[\exp(-\theta Y_t) | Y_0 = y] = \exp\left(-\frac{y\theta}{(t\theta^\beta + 1)^{1/\beta}}\right)(1 + t\theta^\beta)^{-(\beta+1)/\beta}$$

(see Section 6). Moreover, it is possible to start the process $Y$ from the initial state $Y_0 = 0$, and the formula (1.2) continues to hold for $y = 0$. The super-process generalization of this construction was considered for $\beta = 1$ in [11, 14, 15] and for general $\beta$ in [16].

For $\beta = 1$, the conditioned process $Y$ can be described informally as a single "immortal particle" constantly throwing off infinitesimally small masses with each mass then evolving according to the dynamics of the unconditioned process. These infinitesimal masses can be interpreted as the single progenitors of families whose lineage splits from the immortal particle at the birth time of the progenitor and are eventually doomed to extinction. Most such families die immediately, but a rare few live for a noninfinitesimal amount of time. More formally, there is a $\sigma$-finite measure $\nu$ on the space of continuous positive excursion paths

$$\mathcal{E}^0 := \{u \in C(\mathbb{R}_+, \mathbb{R}_+) : u_0 = 0 \ \& \ \exists \gamma > 0 \text{ s.t. } u_t > 0 \Leftrightarrow 0 < t < \gamma\},$$

such that if $\Pi$ is a Poisson point process on $\mathbb{R}_+ \times \mathcal{E}^0$ with intensity $\lambda \otimes \nu$ where $\lambda$ is Lebesgue measure, and $(\bar{X}_t)_{t \geq 0}$ is an independent copy of $X$ begun at $\bar{X}_0 = y$, then the process

$$(1.3) \qquad \left( \bar{X}_t + \sum_{(s,u) \in \Pi} u_{(t-s) \vee 0} \right)_{t \geq 0}$$

has the same distribution as $(Y_t)_{t \geq 0}$ begun at $Y_0 = y$ (see [15]). A point $(s, u) \in \Pi$ corresponds to a family that grows to nonnegligible size; the time $s$ records the moment the family splits off from the immortal particle, and the value $u_r$ of the trajectory $u$ gives the size of the family $r$ units of time after it split off. The family becomes extinct after the period of time $\gamma(u) := \inf\{r > 0 : u(t) = 0, \forall t > r\}$. The $\sigma$-finite measure $\nu$ may be identified explicitly, but it suffices to remark here that it is Markovian with transition probabilities the same as those of the unconditioned process $X$—in other words, $\nu$ arises from a family of entrance laws for the semigroup of $X$. The process $(\bar{X}_t)_{t \geq 0}$ records the population mass due to descendants of individuals other than the immortal particle who are present at time 0.

An analogous description of the conditioned process $Y$ for the case $\beta \in (0, 1)$ is presented in [16]. There is again a single immortal lineage, but now families split off from that lineage with a noninfinitesimal initial size, reflecting the heavy-tailed offspring distributions underlying these models. More precisely, a decomposition similar to (1.3) holds, but the Poisson point process $\Pi$ is now on $\mathbb{R}_+ \times \mathcal{E}$ where $\mathcal{E} = \{u \in D(\mathbb{R}_+, \mathbb{R}_+) : \exists \gamma > 0 \text{ s.t. } u_t > 0 \Leftrightarrow 0 \leq t < \gamma\}$, the set of càdlàg paths starting above zero that eventually hit zero. The nondecreasing process $(M_t)_{t \geq 0}$ where

$$M_t := \sum_{(s,u) \in \Pi \cap [0,t] \times \mathcal{E}} u_0$$

is the total of the initial family sizes that split off from the immortal particle in the time interval $[0, t]$. It is a stable subordinator of index $\beta$.

Suppose now that $\beta \in (0, 1]$ is arbitrary. Take $Y_0 = 0$ so that $\bar{X}_t \equiv 0$ in the decomposition (1.3) and all "individuals" belong to families that split off from the immortal particle at times $s \geq 0$. Extend the definition of $\gamma(u)$ given above for $u \in \mathcal{E}_0$ to $u \in \mathcal{E}$ in the obvious way. The individuals, besides the immortal particle alive at time $t > 0$, belong to families that correspond to the subset $\mathcal{A}_t := \{(s,u) \in \Pi : 0 \leq t - s < \gamma(u)\}$ of the random set $\Pi$. At time $t$, the amount of time since the most recent common ancestor of the entire population lived is $A_t := \sup\{t - s : (s,u) \in \mathcal{A}_t\}$. As depicted in Figure 1, the MRCA age process $(A_t)_{t \geq 0}$ has saw-tooth sample paths that drift up with slope 1 until the current oldest family is extinguished, at which time they jump downward to the age of the next-oldest family.

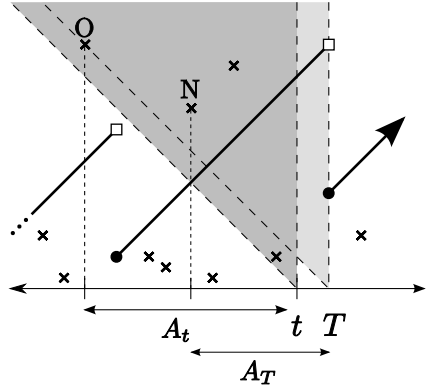

FIG. 1. *The points of $\Lambda$ are marked "x"; the sample path of the process $A$ is drawn with a solid line; the oldest extant family at time $t$ is represented by the point in $\Lambda$ labeled "O" and the left-leaning wedge $\triangle(t,0)$ is the darkly shaded region with apex at $(t,0)$. At time $T$, the family represented by $O$ will die out and the family represented by $N$ will contain the MRCA because $O$ is on the boundary of (the lightly shaded) $\triangle(T,0)$. If the coordinates of $O$ are $(u,y)$ and the coordinates of $N$ are $(v,z)$, then the size of the jump that $A$ makes at $T$ is $A_T - A_{T-} = u - v = (T-v) - (T-u) = (T-v) - y$.*

It is not necessary to know the Poisson point process $\Pi$ in order to construct the MRCA age process $(A_t)_{t \geq 0}$. Clearly, it is enough to know the point process $\Lambda$ on $\mathbb{R}_+ \times \mathbb{R}_{++}$ given by $\Lambda := \{(s, \gamma(u)) : (s,u) \in \Pi\}$. Indeed, if we define the *left-leaning wedge with apex at* $(t,x)$ to be the set

$$(1.4) \qquad \triangle(t,x) := \{(u,v) \in \mathbb{R}^2 : u < t \ \& \ u + v > t + x\},$$

then

$$A_t = t - \inf\{s : \exists x > 0 \text{ s.t. } (s,x) \in \Lambda \cap \triangle(t,0)\}$$

(see Figure 1).

Note that $\Lambda$ is a Poisson point process with intensity $\lambda \otimes \mu$ where $\mu$ is the push-forward of $\nu$ by $\gamma$; that is, $\mu((t,\infty)) = \nu(\{u : \gamma(u) > t\})$. We will show in Section 5 that $\mu((t,\infty)) = (1+\beta)/(\beta t)$.

With these observations in mind, we see that if $\Lambda$ is now *any* Poisson point process on $\mathbb{R}_+ \times \mathbb{R}_{++}$ with intensity $\lambda \otimes \mu$ where $\mu$ is *any* measure on $\mathbb{R}_{++}$ with $\mu(\mathbb{R}_{++}) = \infty$ and $0 < \mu((x,\infty)) < \infty$ for all $x > 0$, then the construction that built $(A_t)_{t \geq 0}$ from the particular point process $\Lambda$ considered above will still apply and produce an $\mathbb{R}_+$-valued process with saw-tooth sample paths. We are therefore led to the following general definition.

DEFINITION 1.1. Let $\Lambda$ be a Poisson point process on $\mathbb{R}_+ \times \mathbb{R}_{++}$ with intensity measure $\lambda \otimes \mu$ where $\lambda$ is Lebesgue measure and $\mu$ is a $\sigma$-finite measure on $\mathbb{R}_{++}$ with $\mu(\mathbb{R}_{++}) = \infty$ and $\mu((x,\infty)) < \infty$ for all $x > 0$. Define $(A_t)_{t \in \mathbb{R}_+}$ by

$$A_t := t - \inf\{s \geq 0 : \exists x > 0 \text{ s.t. } (s,x) \in \Lambda \cap \triangle(t,0)\},$$

where $\triangle(t,0)$ is defined by (1.4), and $A_t = 0$ if $\Lambda \cap \triangle(t,0)$ is empty.

We will suppose from now on that we are in this general situation unless we specify otherwise. We will continue to use terminology appropriate for the genealogical setting and refer to $(A_t)_{t \geq 0}$ as the *MRCA age process* and $\mu$ as the *lifetime measure*. We will assume for convenience that the measure $\mu$ is absolutely continuous with a density $m$ with respect to Lebesgue measure that is positive Lebesgue almost everywhere. It is straightforward to remove these assumptions.

The strong Markov property of the Poisson point processes $\Lambda$ implies that $(A_t)_{t \geq 0}$ is a time-homogeneous strong Markov process. In particular, there is a family of probability distributions $(\mathbb{P}^x)_{x \in \mathbb{R}_+}$ on the space of $\mathbb{R}_+$-valued càdlàg paths with $\mathbb{P}^x$ interpreted in the usual way as the "distribution of $(A_t)_{t \geq 0}$ started from $A_0 = x$." More concretely, the probability measure $\mathbb{P}^x$ is the distribution of the process $(A_t^x)_{t \geq 0}$ defined as follows. Let $\Lambda^x$ be a point process on $[-x, \infty) \times \mathbb{R}_{++}$ that has the distribution of the random point set $\{(t - x, y) : (t, y) \in \Lambda\} \cup \{(-x, Z)\}$ where $Z$ is an independent random variable that is defined on the same probability space $(\Omega, \mathcal{F}, \mathbb{P})$ as the point process $\Lambda$, takes values in the interval $(x, \infty)$ and has distribution

$$\mathbb{P}\{Z \leq z\} = \mu((x, z])/\mu((x, \infty)).$$

Then we can define

$$A_t^x := t - \inf\{s \geq -x : \exists y > 0 \text{ s.t. } (s, y) \in \Lambda^x \cap \triangle(t, 0)\}, \qquad t \geq 0,$$

where we adopt the convention that $A_t^x := 0$ if the set above is empty [see the proof of part (a) of Theorem 1.1 below for further information]. From now on, when we speak of the process $(A_t)_{t \geq 0}$ we will be referring either to the process constructed as in Definition 1.1 from the Poisson process $\Lambda$ defined on some abstract probability space $(\Omega, \mathcal{F}, \mathbb{P})$ or to the canonical process on the space of càdlàg $\mathbb{R}_+$-valued paths equipped with the family of probability measures $(\mathbb{P}^x)_{x \geq 0}$. This should cause no confusion.

We note in passing that the analogue of $(A_t)_{t \geq 0}$ in the constant population size Wright–Fisher setting is not Markov (see Remark 4.1.3 of [26]).

We prove the following properties of the process $(A_t)_{t \geq 0}$ in Section 2.

THEOREM 1.1. (a) *The transition probabilities of the time-homogeneous Markov process $(A_t)_{t \geq 0}$ have an absolutely continuous part*

$$\mathbb{P}^x\{A_t \in dy\} = \frac{\mu((x, x+t])}{\mu((x, \infty))} \exp\left(-\int_y^{x+t} \mu((z, \infty))\, dz\right) \mu((y, \infty))\, dy,$$

*for $y < x + t$ and a single atom*

$$\mathbb{P}^x\{A_t = x + t\} = \frac{\mu((x+t, \infty))}{\mu((x, \infty))}.$$

(b) *The total rate at which the process $(A_t)_{t \geq 0}$ jumps from state $x > 0$ is*

$$\frac{m(x)}{\mu((x, \infty))},$$

*and when the process jumps from state $x > 0$, the distribution of the state to which it jumps is absolutely continuous with density*

$$y \mapsto \exp\left(-\int_y^x \mu((z, \infty))\, dz\right) \mu((y, \infty)), \qquad 0 < y < x.$$

(c) *The probability $\mathbb{P}^0\{\exists t > 0 : A_t = 0\}$ that the process $(A_t)_{t \geq 0}$ returns to the state zero is positive if and only if*

$$\int_0^1 \exp\left(\int_x^1 \mu((y, \infty))\, dy\right) dx < \infty.$$

(d) *If*

$$\int_1^\infty \exp\left(-\int_1^x \mu((y, \infty))\, dy\right) dx = \infty,$$

*then for each $x > 0$ the set $\{t \geq 0 : A_t = x\}$ is $\mathbb{P}^x$-almost surely unbounded. Otherwise, $\lim_{t \to \infty} A_t = \infty$, $\mathbb{P}^x$-almost surely, for all $x \geq 0$.*

(e) *A stationary distribution $\pi$ exists for the process $(A_t)_{t \geq 0}$ if and only if*

$$\int_1^\infty \mu((z, \infty))\, dz < \infty,$$

*in which case it is unique, and*

$$\pi(dx) = \mu((x, \infty)) \exp\left(-\int_x^\infty \mu((z, \infty))\, dz\right) dx.$$

(f) *If $(A_t)_{t \geq 0}$ has a stationary distribution $\pi$, then*

$$d_{\mathrm{TV}}(\mathbb{P}^x\{A_t \in \cdot\}, \pi) \leq 1 - \exp\left(-\int_{t+x}^\infty \mu((y, \infty))\, dy\right) \times \frac{\mu([x, x+t))}{\mu([x, \infty))},$$

*where $d_{\mathrm{TV}}$ denotes the total variation distance. In particular, the distribution of $A_t$ under $\mathbb{P}^x$ converges to $\pi$ in total variation as $t \to \infty$.*

Specializing Theorem 1.1 to the case when $A$ is the MRCA age process of the conditioned critical $(1+\beta)$-stable continuous state branching process gives parts (a) to (d) of the following result. Part (e) follows from an observation that a space–time rescaling of this MRCA age process is a time-homogeneous Markov process that arises from another Poisson process by the general MRCA age construction of Definition 1.1. The proof is in Section 5.

COROLLARY 1.1. *Suppose that $A$ is the MRCA age process associated with the critical $(1+\beta)$-stable continuous state branching process.*

(a) *The transition probabilities of the process $A$ have an absolutely continuous part*

$$\mathbb{P}^x\{A_t \in dy\} = \frac{(1+\beta)ty^{1/\beta}}{\beta(x+t)^{2+1/\beta}}\,dy, \qquad 0 < y < x+t,$$

*and a single atom*

$$\mathbb{P}^x\{A_t = x+t\} = \frac{x}{x+t}.$$

(b) *The total rate at which the process $A$ jumps from the state $x>0$ is $1/x$, and when it jumps from state $x>0$, the distribution of the state to which it jumps is absolutely continuous with density*

$$(1+1/\beta)\frac{y^{1/\beta}}{x^{1+1/\beta}}, \qquad 0 < y < x.$$

(c) *The probability $\mathbb{P}^0\{\exists t>0: A_t = 0\}$ that the process $A$ returns to the state zero is $0$.*

(d) *For each $x \geq 0$, $\lim_{t\to\infty} A_t = \infty$, $\mathbb{P}^x$-almost surely.*

(e) *The process*

$$(e^{-t}A_{e^t})_{t \in \mathbb{R}}$$

*indexed by the whole real line is a time-homogeneous Markov process under $\mathbb{P}^x$ for any $x \geq 0$, and it is stationary when $x=0$. Moreover, $A_t/t$ converges in distribution to the $\mathrm{Beta}(1+1/\beta, 1)$ distribution as $t \to \infty$ under $\mathbb{P}^x$ for any $x \geq 0$, and $A_t/t$ has the $\mathrm{Beta}(1+1/\beta, 1)$ distribution for all $t > 0$ when $x = 0$.*

Note that the sample paths of $(A_t)_{t \geq 0}$ have local "peaks" immediately before jumps and local "troughs" immediately after. We investigate the discrete time Markov chain of successive pairs of peaks and troughs in Section 4. We also consider the jump heights and inter-jump intervals and describe an interesting duality between these sequences in Section 3.

Finally, recall that the Bessel-squared process in dimension $\gamma$, where $\gamma$ is an arbitrary nonnegative real number, is the $\mathbb{R}_+$-valued diffusion process with infinitesimal generator $2x\,d^2/dx^2 + \gamma\,d/dx$. When $\gamma$ is a positive integer, such a process has the same distribution as the square of the Euclidean norm of a Brownian motion in $\mathbb{R}^\gamma$. Feller's critical continuous state branching process is thus the zero-dimensional Bessel-squared process, modulo a choice of scale in time or space. It was shown in Example 3.5 of [28] that for $0 \le \gamma < 2$, the Bessel-squared process with dimension $\gamma$ conditioned on never hitting zero is the Bessel-squared process with dimension $4 - \gamma$. Thus for $\beta = 1$, the conditioned process $Y$ is the four-dimensional Bessel-squared process. We introduce a new family of processes in Section 6 that are also indexed by a nonnegative real parameter and play the role of the Bessel-squared family for values of $\beta$ other than 1. These processes will be studied further in a forthcoming paper.

We end this introduction by commenting on the connections with previous work. First, we may think of each point $(s, x) \in \Lambda$ as a "job" that enters a queue with infinitely many servers at time $s$ and requires an amount of time $x$ to complete. We thus have a classical $M/G/\infty$ queue [31], except we are assuming that the total arrival rate of jobs of all kinds is infinite. With this interpretation, the quantity $A_t$ is how long the oldest job at time $t$ has been in the queue. Properties of such $M/G/\infty$ queues with infinite arrival rates have been studied (see, for example, [10]), but the age of the oldest job does not appear to have been studied in this context.

Second, note that the process $(A_t)_{t \ge 0}$ is an example of a *piecewise deterministic Markov process*: it consists of deterministic flows punctuated by random jumps. Such processes were introduced in [5] and studied further in [6] (see also [20], where the nomenclature *jumping Markov processes* is used). The general properties of such processes have been studied further in, for example, [3, 4, 7].

Last, piecewise deterministic Markov processes like $(A_t)_{t \ge 0}$ that have periods of linear increase interspersed with random jumps have been used to model many phenomena, such as stress in an earthquake zone [2], congestion in a data transmission network [8] and growth-collapse [1]. They also have appeared in the study of the additive coalescent [12] and $\mathbb{R}$-tree-valued Markov processes [13].

**2. Proof of Theorem 1.1.** (a) Suppose that $(A_t)_{t \ge 0}$ is constructed from $\Lambda$ as in Definition 1.1. For $s \ge x$ consider the conditional distribution of $A_{s+t}$ given $A_s = x$. The condition $A_s = x$ is equivalent to the requirements that there is a point $(s-x, Z)$ in $\Lambda$ for some $Z > x$ and that, furthermore, there are no points of $\Lambda$ in the left-leaning wedge $\triangle(s-x, x)$. The conditional probability of the event $\{Z > z\}$, given $A_s = x$, is $\mu((z, \infty))/\mu((x, \infty))$ for $z \geq x$. If $Z > x + t$, then $A_{s+t} = x + t$. Otherwise, $A_{s+t} < x + t$. The second claim of part (a) follows immediately.

Now consider $\mathbb{P}\{A_{s+t} \in dy | A_s = x\}$ for $y < x + t$. This case is depicted in Figure 2. By construction, $A_{s+t} = y$ if and only if there is a point $(s+t-y, W) \in \Lambda$ for some $W > y$, and there are no points of $\Lambda$ in $\triangle(s+t-y, y)$. From above, the condition $A_s = x$ requires there to be no points of $\Lambda$ in the wedge $\triangle(s-x, x)$ (the lightly shaded region in Figure 2). Therefore, if $A_s = x$, then $A_{s+t} = y$ if and only if $Z \leq x + t$, there are no points of $\Lambda$ in the darkly shaded region of Figure 2 and there is a point of $\Lambda$ of the form $(s+t-y, w)$ with $w > y$.

Now, the conditional probability of the even $\{Z \leq x + t\}$ given $A_s = x$ is

$$\frac{\mu((x, x+t])}{\mu((x, \infty))}.$$

The probability that no points of $\Lambda$ are in the darkly shaded region is

$$\exp\left(-\int_y^{x+t} \mu((u, \infty)) \, du\right).$$

The probability that $\Lambda$ has a point in the infinitesimal region $[s+t-y, s+t-y+dy] \times (y, \infty)$ is $\mu((y, \infty)) \, dy$. Multiplying these three probabilities together gives the first claim of part (a).

(b) Both claims follow readily by differentiating the formulae in (a) at $t = 0$. It is also possible to argue directly from the representation in terms of $\Lambda$.

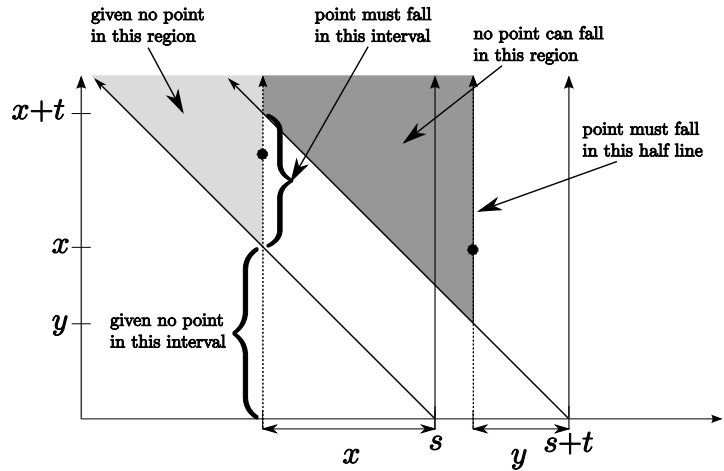

FIG. 2.    The computation of $\mathbb{P}\{A_{s+t} \in dy | A_s = x\}$.

(c) Suppose that $(A_t)_{t\geq 0}$ is constructed from $\Lambda$ as in Definition 1.1. Note that

$$\{t > 0 : A_t = 0\} = \mathbb{R}_{++} \setminus \bigcup_{(t,x)\in\Lambda} (t, t+x).$$

The question of when such a "Poisson cut out" random set is almost surely empty was asked in [23], and the necessary and sufficient condition presented in part (c) is simply the one found in [30] (see also [18]).

(d) For $x > 0$, the random set $\{t \geq 0 : A_t = x\}$ is a discrete regenerative set under $\mathbb{P}^x$ that is not almost surely equal to $\{0\}$. Hence this set is a renewal process with an inter-arrival time distribution that possibly places some mass at infinity; in which case, the number of arrivals is almost surely finite with a geometric distribution, and the set is almost surely bounded.

Suppose that the set $\{t \geq 0 : A_t = x\}$ is $\mathbb{P}^x$-almost surely unbounded for some $x > 0$. Let $0 = T_0 < T_1 < \cdots$ be the successive visits to $x$. It is clear that $\mathbb{P}^x\{\exists t \in [0, T_1] : A_t = y\} > 0$ for any choice of $y > 0$ and hence, by the strong Markov property, the set $\{t \geq 0 : A_t = y\}$ is also $\mathbb{P}^x$-almost surely unbounded. Another application of the strong Markov property establishes that the set $\{t \geq 0 : A_t = y\}$ is $\mathbb{P}^y$-almost surely unbounded. Thus the set $\{t \geq 0 : A_t = x\}$ is either unbounded $\mathbb{P}^x$-almost surely for all $x > 0$ or bounded $\mathbb{P}^x$-almost surely for all $x > 0$.

Recall that away from the set $\{t \geq 0 : A_t = 0\}$ the sample paths of $A$ are piecewise linear with slope 1. It follows from the coarea formula (see, e.g., Section 3.8 of [24]) that

$$\int_0^\infty f(A_t)\,dt = \int_0^\infty f(y) \#\{t > 0 : A_t = y\}\,dy$$

for a Borel function $f : \mathbb{R}_+ \to \mathbb{R}_+$. Hence by Fubini's theorem,

$$\int_0^\infty f(y)\mathbb{E}^x[\#\{t > 0 : A_t = y\}]\,dy = \int_0^\infty f(y) \int_0^\infty \frac{\mathbb{P}^x\{A_t \in dy\}}{dy}\,dt\,dy$$

for any $x > 0$. It follows from the continuity of the transition probability densities that

$$\mathbb{E}^x[\#\{t > 0 : A_t = y\}] = \int_0^\infty \frac{\mathbb{P}^x\{A_t \in dy\}}{dy}\,dt$$

for all $x, y > 0$ and, in particular, that the expected number of returns to $x > 0$ under $\mathbb{P}^x$ is

$$\int_0^\infty \frac{\mathbb{P}^x\{A_t \in dx\}}{dx}\,dt.$$

Using the expression from part (a) and the argument above, this quantity is infinite, and hence the number of visits is $\mathbb{P}^x$-almost surely infinite, if and

only if
$$\int_1^\infty \exp\left(-\int_1^u \mu((y,\infty))\,dy\right)du = \infty.$$

If the set $\{t \geq 0 : A_t = x\}$ is $\mathbb{P}^x$-almost surely bounded for all $x > 0$, then by an argument similar to the above, the set $\{t \geq 0 : A_t = y\}$ is $\mathbb{P}^x$-almost surely bounded for all $x, y > 0$. It follows that, for all $x > 0$, $\mathbb{P}^x$-almost surely all of the sets $\{t \geq 0 : A_t = y\}$ are finite. This implies that $\lim_{t \to \infty} A_t$ exists $\mathbb{P}^x$-almost surely, and the limit takes values in the set $\{0, \infty\}$. However, it is clear from the Poisson process construction that 0 does not occur as a limit with positive probability.

(e) Suppose there exists a probability measure $\pi$ on $\mathbb{R}_+$ such that
$$\int_{\mathbb{R}_+} \mathbb{P}\{A_t \in dy | A_0 = x\}\pi(dx) = \pi(dy), \qquad y \in \mathbb{R}_+.$$

Taking $t \to \infty$ in part (a) gives
$$\pi(dy) = \lim_{t \to \infty} \int_{\mathbb{R}_+} \frac{\mu((x, x+t])}{\mu((x, \infty))}$$
$$\times \exp\left(-\int_y^{x+t} \mu((u, \infty))\,du\right)\pi(dx)\mu((y, \infty))\,dy$$
$$= \begin{cases} 0, & \text{if } \int_y^\infty \mu((u, \infty))\,du = \infty, \\ \exp\left(-\int_y^\infty \mu((u, \infty))\,du\right)\mu((y, \infty))\,dy, & \text{otherwise.} \end{cases}$$

Therefore, a stationary probability distribution exists if and only if $\int_y^\infty \mu((u, \infty))\,du < \infty$, and if a stationary distribution exists, then it is unique.

(f) It will be useful to begin with a concrete construction of a stationary version of the process $A$ in terms of a Poisson point process. Suppose that $\int_x^\infty \mu((u, \infty))\,du < \infty$ for all $x > 0$, so that a stationary distribution exists. Let $\Lambda^\leftrightarrow$ be a Poisson point process on $\mathbb{R} \times \mathbb{R}_{++}$ with intensity measure $\lambda \otimes \mu$. Define $(A_t^\leftrightarrow)_{t \in \mathbb{R}}$ by
$$A_t^\leftrightarrow := t - \inf\{s : \exists x > 0 \text{ s.t. } (s, x) \in \Lambda^\leftrightarrow \cap \triangle(t, 0)\}.$$

The condition on $\mu$ ensures that almost surely any wedge $\triangle(t, x)$ with $x > 0$ will contain only finitely many points of $\Lambda^\leftrightarrow$, and so $(A_t^\leftrightarrow)_{t \in \mathbb{R}}$ is well defined. The process $(A_t^\leftrightarrow)_{t \in \mathbb{R}}$ is stationary and Markovian, with the same transition probabilities as $(A_t)_{t \geq 0}$.

Recall the construction of the process $A^x$ started at $x$ for $x > 0$ that was described preceding the statement of Theorem 1.1. Construct the point process $\Lambda^x$ that appears there by setting $\Lambda^x := \{(t, y) \in \Lambda^\leftrightarrow : t > -x\} \cup \{(-x, Z)\}$ where $Z$ is an independent random variable with values in the interval $(x, \infty)$ and distribution $\mathbb{P}\{Z > z\} = \mu((z, \infty))/\mu((x, \infty))$.

By construction, $A_t^x = A_t^{\leftrightarrow}$ for all $t \geq T$ where $T$ is the death time of all families alive at time $-x$ in either process:

$$T := \inf\{t > 0 : Z \leq t + x \text{ and } \Lambda^{\leftrightarrow} \cap \triangle(-x, t+x) = \varnothing\}.$$

Thus

$$\begin{aligned}
d_{\mathrm{TV}}(\mathbb{P}^x\{A_t \in \cdot\}, \pi) &\leq \mathbb{P}\{A_t^x \neq A_t^{\leftrightarrow}\} \\
&\leq \mathbb{P}\{T > t\} \\
&= 1 - \mathbb{P}\{Z \leq t + x\}\mathbb{P}\{\Lambda^{\leftrightarrow} \cap \triangle(-x, t+x) = \varnothing\}
\end{aligned}$$

and part (f) follows.

**3. Duality and time-reversal.** Suppose in this section that $\int_x^{\infty} \mu((y, \infty)) \, dy < \infty$ for all $x > 0$, so that, by part (e) of Theorem 1.1, the process $A$ has a stationary distribution. Let $(A_t^{\leftrightarrow})_{t \in \mathbb{R}}$ be the stationary Markov process with the transition probabilities of $A$ that was constructed from the Poisson point process $\Lambda^{\leftrightarrow}$ in the proof of part (f) of Theorem 1.1.

Define the *dual process* $(\widehat{A}_t^{\leftrightarrow})_{t \in \mathbb{R}}$ by $\widehat{A}_t^{\leftrightarrow} := \inf\{s > 0 : \triangle(t, s) \cap \Lambda^{\leftrightarrow} = \varnothing\}$. See Figures 3 and 4. Thus $\widehat{A}_t^{\leftrightarrow}$ is the amount of time that must elapse after time $t$ until all families alive at time $t$ have died out, or, equivalently, until the MRCA for the population lived at some time after $t$. The càdlàg $\mathbb{R}_+$-valued process $(\widehat{A}_t^{\leftrightarrow})_{t \in \mathbb{R}}$ has saw-tooth sample paths that drift down with slope $-1$ between upward jumps.

PROPOSITION 3.1. *The dual process, $(\widehat{A}_t^{\leftrightarrow})_{t \in \mathbb{R}}$, has the same distribution as the time-reversed process, $(\bar{A}_t^{\leftrightarrow})_{t \in \mathbb{R}}$, where $\bar{A}_t := \lim_{u \downarrow t} A_{-u}^{\leftrightarrow}$.*

PROOF. Define bijections $\varphi$ and $\sigma$ of $\mathbb{R} \times \mathbb{R}_{++}$ by $\varphi(t, x) := (t + x, x)$ and $\sigma(t, x) := (-t, x)$. In the usual manner, we may also think of $\varphi$ and $\sigma$ as mapping subsets of $\mathbb{R} \times \mathbb{R}_{++}$ to other subsets of $\mathbb{R} \times \mathbb{R}_{++}$. Note that $\varphi$ maps left-leaning wedges to right-leaning wedges, and $\sigma$ maps right-leaning wedges to left-leaning wedges. Thus the composition $\sigma \circ \varphi$ maps left-leaning wedges to left-leaning wedges. More precisely,

$$\begin{aligned}
\sigma \circ \varphi(\triangle(t, x)) &\\
&= \{\sigma \circ \varphi(s, y) : (s, y) \in \triangle(t, x)\} \\
(3.1) \quad &= \{(-(s+y), y) : s < t \ \& \ s + y > t + x\} \\
&= \{(u, v) : u < -(t+x) \ \& \ u + v > -t = -(t+x) + x\} \\
&= \triangle(\sigma \circ \varphi(t, x)),
\end{aligned}$$

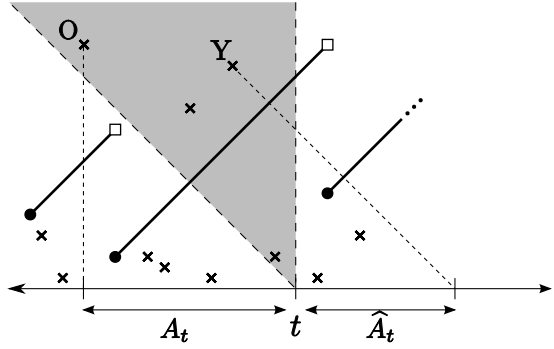

FIG. 3. *The process $A_t^{\leftrightarrow}$ and the dual process $\widehat{A}_t^{\leftrightarrow}$. The points in the shaded area represent the families alive at time $t$, and the solid line is the sample path of $A_t^{\leftrightarrow}$. The point marked "O" is the oldest living family at time $t$; the point marked "Y" is the family extant at time $t$ that will live the longest into the future.*

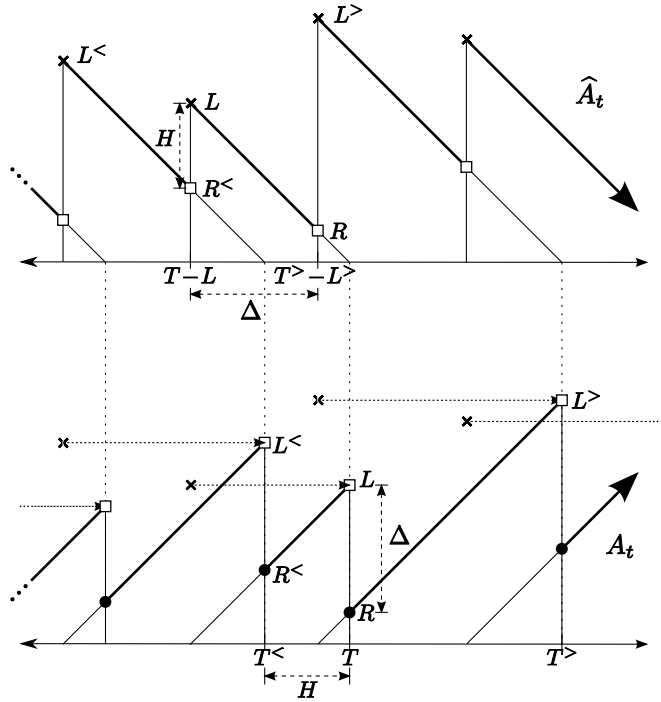

FIG. 4. *The (coupled) processes $A^{\leftrightarrow}$ and $\widehat{A}^{\leftrightarrow}$ for the same set of family lifetimes (the "x"s). The paths of $\widehat{A}^{\leftrightarrow}$ begin at points in $\Gamma(\Lambda^{\leftrightarrow})$ and descend; while the paths of $A^{\leftrightarrow}$ ascend to points in $\varphi(\Gamma(\Lambda^{\leftrightarrow}))$. The mapping $\varphi$ is shown by the horizontal dotted lines in the lower diagram.*

where we stress that $\triangle(\sigma \circ \varphi(t,x))$ is the left-leaning wedge with its apex at the point $\sigma \circ \varphi(t,x) \in \mathbb{R} \times \mathbb{R}_{++}$. Define a map $\Gamma$ that takes a subset of $\mathbb{R} \times \mathbb{R}_{++}$ and returns another such subset by

$$\Gamma(S) := \{(t,x) \in S : \triangle(t,x) \cap S = \varnothing\}.$$

The points of $\Gamma(\Lambda^{\leftrightarrow})$ correspond precisely to those families that at some time will be the oldest surviving family in the population. These points determine the jumps of both the MRCA process and the dual process: as can be seen with the help of Figure 4, the linear segments of the paths of the dual process $(\widehat{A}^{\leftrightarrow}_t)_{t \in \mathbb{R}}$ each begin at a point in $\Gamma(\Lambda^{\leftrightarrow})$ and descend with slope $-1$ whereas the linear segments of the paths of the MRCA process $(A^{\leftrightarrow}_t)_{t \in \mathbb{R}}$ ascend with slope $+1$ to points in $\varphi \circ \Gamma(\Lambda^{\leftrightarrow})$.

This implies that the path of the time-reversed process, $(\bar{A}_t)_{t \in \mathbb{R}}$, begin at points in $\sigma \circ \varphi \circ \Gamma(\Lambda^{\leftrightarrow})$ and descend with slope $-1$, and so the points of $\sigma \circ \varphi \circ \Gamma(\Lambda^{\leftrightarrow})$ determine the path of $(\bar{A}_t)_{t \in \mathbb{R}}$ in the same way that the points of $\Gamma(\Lambda^{\leftrightarrow})$ determine the path of $(\widehat{A}_t)_{t \in \mathbb{R}}$. Therefore, all we need to do is show that $\sigma \circ \varphi \circ \Gamma(\Lambda^{\leftrightarrow})$ has the same distribution as $\Gamma(\Lambda^{\leftrightarrow})$.

To see this, first note that, for an arbitrary subset $S \subset \mathbb{R} \times \mathbb{R}_{++}$,

$$\begin{aligned}
\sigma &\circ \varphi \circ \Gamma(S) \\
&= \{\sigma \circ \varphi(t,x) : (t,x) \in S \text{ and } \triangle(t,x) \cap S = \varnothing\} \\
&= \{(t,x) : \varphi^{-1} \circ \sigma^{-1}(t,x) \in S \text{ and } \triangle(\varphi^{-1} \circ \sigma^{-1}(t,x)) \cap S = \varnothing\} \\
&= \{(t,x) : \varphi^{-1} \circ \sigma^{-1}(t,x) \in S \text{ and } \varphi^{-1} \circ \sigma^{-1}(\triangle(t,x)) \cap S = \varnothing\} \\
&= \{(t,x) : (t,x) \in \sigma \circ \varphi(S) \text{ and } \triangle(t,x) \cap \sigma \circ \varphi(S) = \varnothing\} \\
&= \Gamma \circ \sigma \circ \varphi(S),
\end{aligned}$$

where the third equality follows from identity (3.1). Thus $\sigma \circ \varphi \circ \Gamma(\Lambda^{\leftrightarrow}) = \Gamma \circ \sigma \circ \varphi(\Lambda^{\leftrightarrow})$. However, both maps $\varphi$ and $\sigma$ preserve the measure $\lambda \otimes \mu$ and hence $\sigma \circ \varphi(\Lambda^{\leftrightarrow})$ has the same distribution as $\Lambda^{\leftrightarrow}$. Therefore, $\sigma \circ \varphi \circ \Gamma(\Lambda^{\leftrightarrow})$ has the same distribution as $\Gamma(\Lambda^{\leftrightarrow})$, as required. □

REMARK 3.1. There is an interesting connection between the jump sizes and the inter-jump intervals, stemming from the observation that the paths of $A^{\leftrightarrow}$ and $\widehat{A}^{\leftrightarrow}$ have the same sequences of "trough" and "peak" heights, while the roles of the jump sizes and inter-jump intervals for the two are exchanged. The following situation is depicted in Figure 4. To explain the connection, suppose that $T \in \{t \in \mathbb{R} : A^{\leftrightarrow}_{t-} \neq A^{\leftrightarrow}_t\}$ is a jump time for the process $A^{\leftrightarrow}$. Let $T^< := \sup\{t < T : A^{\leftrightarrow}_{t-} \neq A^{\leftrightarrow}_t\}$ and $T^> := \inf\{t > T : A^{\leftrightarrow}_{t-} \neq A^{\leftrightarrow}_t\}$ be the jump times on either side of $T$. Put $L := A^{\leftrightarrow}_{T-}$ and $R := A^{\leftrightarrow}_T$, and define $L^<$, $R^<$, $L^>$ and $R^>$ as the analogous left limits and values of $A^{\leftrightarrow}$ at the times $T^<$ and $T^>$. Write $\Delta := L - R$ for the size of the jump at time $T$ and $H := T - T^<$ for the length of the time interval since the previous jump. Observe that $T - L$ is a jump time for the dual process $\widehat{A}^{\leftrightarrow}$ with $\widehat{A}^{\leftrightarrow}_{(T-L)-} = R^<$ and $\widehat{A}^{\leftrightarrow}_{(T-L)} = L$. Moreover,

$$H = L - R^< = \widehat{A}^{\leftrightarrow}_{(T-L)} - \widehat{A}^{\leftrightarrow}_{(T-L)-}$$

and

$$\Delta = (T^> - L^>) - (T - L) = \inf\{t > T - L : \widehat{A}^{\leftrightarrow}_{t-} \neq \widehat{A}^{\leftrightarrow}_t\} - (T - L).$$

Note that the map $T \mapsto T - L$ sets up a monotone bijection between the jump times of the process $A^{\leftrightarrow}$ and those of the process $\widehat{A}^{\leftrightarrow}$. It thus follows from Proposition 3.1 that the point processes,

$$\{(T, T - T^<, A^{\leftrightarrow}_{T-} - A^{\leftrightarrow}_T) : A^{\leftrightarrow}_{T-} \neq A^{\leftrightarrow}_T\}$$

and

$$\{(T, A^{\leftrightarrow}_{T-} - A^{\leftrightarrow}_T, T - T^<) : A^{\leftrightarrow}_{T-} \neq A^{\leftrightarrow}_T\},$$

have the same distribution.

**4. Jump chains.** Suppose again that $\int_x^\infty \mu((y, \infty)) \, dy < \infty$ for all $x > 0$ so that the process $A$ has a stationary distribution. Recall the stationary Markov process, $(A^{\leftrightarrow}_t)_{t \in \mathbb{R}}$, with the transition probabilities of $A$ that was constructed from the Poisson point process, $\Lambda^{\leftrightarrow}$, in the proof of part (f) of Theorem 1.1.

For $t \in \mathbb{R}$, denote by $J_t := \inf\{u > 0 : A^{\leftrightarrow}_u \neq A^{\leftrightarrow}_{u-}\}$ the next jump time of $A^{\leftrightarrow}$ after time $t$. Define an increasing sequence of random times, $0 < T_0 < T_1 < \cdots$, by $T_0 := J_0$ and $T_{n+1} := J_{T_n}$ for $n \geq 0$. Put $L_n := A^{\leftrightarrow}_{T_n-}$ and $R_n := A^{\leftrightarrow}_{T_n}$. Thus the sequences $(L_n)_{n=0}^\infty$ and $(R_n)_{n=0}^\infty$ record, respectively, the "peaks" and the "troughs" of the path of $A^{\leftrightarrow}$ that occur between the times $0$ and $\sup_n T_n$.

The next result can be proved along the same lines as part (a) of Theorem 1.1, and we leave the proof to the reader.

PROPOSITION 4.1. *The sequence, $(L_0, R_0, L_1, R_1, \ldots)$, is Markovian with the following transition probabilities:*

$$\mathbb{P}\{R_n \in dy | L_n = x\} = \mu((y, \infty)) \exp\left(-\int_y^x \mu((u, \infty)) \, du\right) dy, \qquad 0 < y \leq x,$$

*and*

$$\mathbb{P}\{L_{n+1} \in dz | R_n = y\} = \frac{m(z)}{\mu((y, \infty))} dz, \qquad z > y.$$

*In particular, the sequence of pairs $((L_n, R_n))_{n=0}^\infty$ is a time-homogeneous Markov chain.*

Now we may compute the transition probabilities of the peak and trough chains. By Proposition 4.1,

$$
\begin{aligned}
\mathbb{P}\{L_{n+1} \in dz | L_n = x\}/dz \\
(4.1) \quad &= \int_0^x \exp\left(-\int_y^x \mu((u,\infty))\,du\right)\mu((y,\infty))\frac{m(z)\mathbf{1}_{y\leq z}}{\mu((y,\infty))}\,dy \\
&= m(z)\int_0^{x\wedge z} \exp\left(-\int_y^x \mu((u,\infty))\,du\right) dy
\end{aligned}
$$

and

$$
\begin{aligned}
\mathbb{P}\{R_{n+1} \in dz | R_n = x\}/dz \\
(4.2) \quad &= \int_x^\infty \frac{m(y)}{\mu((x,\infty))}\mu((z,\infty))\exp\left(-\int_z^y \mu((u,\infty))\,du\right)\mathbf{1}_{y>z}\,dy \\
&= \mu((z,\infty))\int_{x\vee z}^\infty \frac{m(y)}{\mu((x,\infty))}\exp\left(-\int_z^y \mu((u,\infty))\,du\right) dy.
\end{aligned}
$$

It follows from (4.1) that the peak chain, $(L_n)_{n=0}^\infty$, is $\lambda$-*irreducible* where $\lambda$ is Lebesgue measure on $\mathbb{R}_{++}$. That is, if $A$ is a Borel subset of $\mathbb{R}_{++}$ with $\lambda(A) > 0$, then, for any $x \in \mathbb{R}_{++}$, there is positive probability that the the peak chain begun at $x$ will hit $A$ at some positive time (see Chapter 4 of [25] for more about this notion of irreducibility). It follows that the peak chain is either *recurrent*, in the sense that

$$\sum_{n=0}^\infty \mathbb{P}\{L_n \in A | L_0 = x\} = \infty$$

for all $x \in \mathbb{R}_{++}$ and all Borel subsets of $A \subseteq \mathbb{R}_{++}$ with $\lambda(A) > 0$, or it is *transient*, in the sense that there is a countable collection of Borel sets $(A_j)_{j=1}^\infty$ and finite constants $(M_j)_{j=1}^\infty$ such that $\bigcup_{j=1}^\infty A_j = \mathbb{R}_{++}$ and

$$\sup_{x\in\mathbb{R}_{++}} \sum_{n=0}^\infty \mathbb{P}\{L_n \in A_j | L_0 = x\} \leq M_j$$

(see Theorem 8.0.1 of [25]).

The peak chain is *strong Feller*; that is, the function,

$$x \mapsto \mathbb{E}[f(L_{n+1})|L_n = x],$$

is continuous for any bounded Borel function $f$. Also, because the support of $\lambda$ is all of $\mathbb{R}_{++}$, if the peak chain is recurrent, then each point $x$ of $\mathbb{R}_{++}$ is *topologically recurrent* in the sense that

$$\sum_{n=0}^\infty \mathbb{P}\{L_n \in U | L_0 = x\} = \infty$$

for every open neighborhood $U$ of $x$. Hence, by Theorem 9.3.6 of [25], if the peak chain is recurrent, then it is *Harris recurrent*, which means that given any Borel set $A$ with $\lambda(A) > 0$, the chain visits $A$ infinitely often almost surely starting from any $x$. Moreover, the chain is recurrent (equivalently, Harris recurrent) if and only if it is *nonevanescent*; that is, started from any $x$, there is zero probability that the chain converges to 0 or $\infty$ (see Theorem 9.2.2 of [25]).

If the peak chain is recurrent (equivalently, Harris recurrent or nonevanescent), then it has an invariant measure that is unique up to constant multiples (see Theorem 10.4.4 of [25]). If the invariant measure has finite mass, so that it can be normalized to be a probability measure, then the chain is said to be *positive*, otherwise the chain is said to be *null*.

Conversely, if the peak chain has an invariant probability measure, then it is recurrent (equivalently, Harris recurrent or nonevanescent) (see Proposition 10.1.1 of [25]).

All of the remarks we have just made for the peak chain apply equally to the trough chain $(R_n)_{n=0}^\infty$. Recall that we are in the situation when $A$ has a stationary version, so the transience or recurrence of $L$ and $R$ depends on their behavior near zero.

PROPOSITION 4.2. *Consider the two Markov chains, $(L_n)_{n=0}^\infty$ and $(R_n)_{n=0}^\infty$.*

(a) *Both chains are transient if and only if*
$$\int_0^1 \exp\left(\int_x^1 \mu((y,\infty))\,dy\right) dx < \infty.$$

(b) *Both chains are positive recurrent if and only if*
$$\int_0^1 m(x) \exp\left(-\int_x^1 \mu((y,\infty))\,dy\right) dx < \infty.$$

(c) *Both chains are null recurrent if and only if both*
$$\int_0^1 \exp\left(\int_x^1 \mu((y,\infty))\,dy\right) dx = \infty$$
*and*
$$\int_0^1 m(x) \exp\left(-\int_x^1 \mu((y,\infty))\,dy\right) dx = \infty.$$

PROOF. Consider the set $\mathcal{Z} := \{t \in \mathbb{R} : A_t^{\leftrightarrow} = 0\}$. It follows from part (c) of Theorem 1.1 that $\mathbb{P}\{\mathcal{Z} \neq \varnothing\} > 0$ if and only if

(4.3) $$\int_0^1 \exp\left(\int_x^1 \mu(y,\infty)\,dy\right) dx < \infty.$$

By the stationarity of $(A_t^{\leftrightarrow})_{t \in \mathbb{R}}$ and the nature of its sample paths, it is clear that for $x > 0$, the set $\{t \in \mathbb{R} : A_t^{\leftrightarrow} = x\}$ is unbounded above and below almost surely [this also follows from part (d) of Theorem 1.1]. It follows from a simple renewal argument that if (4.3) holds, then $\mathcal{Z}$ is unbounded above and below almost surely.

Because the paths of $(A_t^{\leftrightarrow})_{t \in \mathbb{R}}$ increase with slope 1 in the intervals $[T_n, T_{n+1})$, it follows that if condition (4.3) holds, then $\lim_{n \to \infty} T_n = \inf\{t > 0 : A_t^{\leftrightarrow} = 0\} < \infty$ almost surely, and $\lim_{n \to \infty} L_n = \lim_{n \to \infty} R_n = 0$ almost surely. In this case, both chains are evanescent, and hence transient.

On the other hand, if (4.3) does not hold, then $\lim_{n \to \infty} T_n = \infty$. Moreover, the set $\{t \in \mathbb{R} : A_t^{\leftrightarrow} = x\}$ is almost surely unbounded above and below for any $x > 0$, as we observed above. If we split the path of $(A_t^{\leftrightarrow})_{t \in \mathbb{R}}$ into excursions away from $x$, then each excursion interval will contain only finitely many jumps almost surely and, because the excursions are independent and identically distributed, it cannot be the case that $L_n$ or $R_n$ converges to 0 or $\infty$ with positive probability. Thus, both chains are nonevanescent and hence recurrent.

It is clear from (4.1) that the kernel giving the transition densities of the peak chain $(L_n)_{n=0}^{\infty}$ is self-adjoint with respect to the measure having density

$$p(x) = m(x) \exp\left(-\int_x^{\infty} \mu((u, \infty)) \, du\right)$$

with respect to Lebesgue measure, and so this measure is invariant for the peak chain. Clearly, $\int_0^{\infty} p(x) \, dx < \infty$ if and only if the condition in part (b) holds in which case the peak chain is positive recurrent. Otherwise, the peak chain is either null recurrent or transient, and so part (a) shows that the peak chain is null recurrent if the two conditions in part (c) hold.

Similarly, it is clear from (4.2) that the kernel giving the transition densities of the trough chain $(R_n)_{n=0}^{\infty}$ is self-adjoint with respect to the measure having density

$$q(x) = \mu((x, \infty))^2 \exp\left(-\int_x^{\infty} \mu((u, \infty)) \, du\right)$$

with respect to Lebesgue measure, and so this measure is invariant for the trough chain. An integration by parts shows that $\int_0^{\infty} q(x) \, dx < \infty$ if and only if the condition in part (b) holds, and so the trough chain is positive if and only if the peak chain is positive. Alternatively, we can simply observe from Proposition 4.1 that integrating the conditional probability kernel of $R_n$ given $L_n$ against an invariant probability measure for the peak chain gives an invariant measure for the trough chain, and integrating the conditional probability kernel of $L_{n+1}$ given $R_n$ against an invariant probability measure for the trough chain chain gives an invariant measure for the trough chain, so that one chain is positive recurrent if and only if the other is.  □

REMARK 4.1. If $m(x) = \alpha x^{-2}$ for $x \in (0,1]$, then both the peak and trough chains are:

(1) transient $\Leftrightarrow 0 < \alpha < 1$;
(2) null recurrent $\Leftrightarrow \alpha = 1$;
(3) positive recurrent $\Leftrightarrow \alpha > 1$.

REMARK 4.2. It follows from parts (b) and (e) of Theorem 1.1 that the stationary point process $\{t \in \mathbb{R} : A_{t-}^{\leftrightarrow} \neq A_t^{\leftrightarrow}\}$ has intensity

$$\rho := \int m(x) \exp\left(-\int_x^\infty \mu((u,\infty))\,du\right) dx$$

and so the peak and trough chains are positive recurrent if and only if $\rho$ is finite. Suppose that $\rho$ is finite and consider the point process

$$\Xi := \{(t, A_{t-}^{\leftrightarrow}, A_t^{\leftrightarrow}) \in \mathbb{R} \times \mathbb{R}_+ \times \mathbb{R}_+ : A_{t-}^{\leftrightarrow} \neq A_t^{\leftrightarrow}\}.$$

The companion *Palm* point process $\Upsilon$ has its distribution defined by

$$\mathbb{P}\{\Upsilon \in \cdot\} = \rho^{-1} \mathbb{E}\left[\sum_{\{n\,:\,0 \leq T_n \leq 1\}} \mathbf{1}\{\theta_{T_n} \Xi \in \cdot\}\right],$$

where $\theta_s B = \{(t-s, \ell, r) : (t, \ell, r) \in B\}$ for $B \subset \mathbb{R} \times \mathbb{R}_+ \times \mathbb{R}_+$. Enumerate the points of $\Upsilon$ as $((\tilde{T}_n, \tilde{L}_n, \tilde{R}_n))_{n \in \mathbb{Z}}$ where $\cdots < \tilde{T}_{-1} < \tilde{T}_0 = 0 < \tilde{T}_1 < \cdots$. A fundamental result of Palm theory for stationary point processes says that the random sequence $((\tilde{T}_n - \tilde{T}_{n-1}, \tilde{L}_n, \tilde{R}_n))_{n \in \mathbb{Z}}$ is stationary and that the distribution of the point process $\Xi$ may be reconstructed from the distribution of this sequence (see, e.g., Theorem 12.3.II of [9] or [21]). It is clear that the stationary random sequences $(\tilde{L}_n)_{n=0}^\infty$ and $(\tilde{R}_n)_{n=0}^\infty$ have the same distribution as the peak and trough chains started in their respective stationary distributions.

**5. The $(1+\beta)$-stable MRCA process.** In this section we specialize to the motivating example of the MRCA process of a critical $(1+\beta)$-stable continuous state branching process conditioned to live forever. Recall that the unconditioned continuous state branching process has Laplace transforms (1.1), and the conditioned process has Laplace transforms (1.2). For $\beta = 1$, the unconditioned process has generator $x\frac{\partial^2}{\partial x^2}$ with this choice of time scale.

As mentioned in the Introduction, the set of points $(t,x) \in \mathbb{R}_+ \times \mathbb{R}_{++}$, where $t$ is the time that a family splits from the immortal lineage and $x$ is its total lifetime, is a Poisson point process with intensity measure $\lambda \otimes \mu$ for some $\sigma$-finite measure $\mu$.

LEMMA 5.1. *The lifetime measure $\mu$ associated with the critical $(1+\beta)$-branching process conditioned on nonextinction is given by*

$$\mu((x,\infty)) = \frac{1+\beta}{\beta x}, \qquad x > 0.$$

PROOF. As we remarked in the Introduction,

(5.1) $$\mathbb{P}\{X_t > 0 | X_0 = x\} = 1 - \exp\left(-\frac{x}{t^{1/\beta}}\right).$$

First consider the case of $\beta = 1$. Recall from the Introduction that if $\Pi$ is a Poisson point process on $\mathbb{R}_+ \times \mathcal{E}^0$ with intensity $\lambda \otimes \nu$, then

(5.2) $$\left(\sum_{(s,u)\in\Pi} u_{(t-s)\vee 0}\right)_{t\geq 0}$$

has the same distribution as the conditioned process, $(Y_t)_{t\geq 0}$, with $Y_0 = 0$, and recall that $\mu$ is the push-forward of $\nu$ by the total lifetime function $\gamma$. Also,

(5.3) $$\left(\sum_{(s,u)\in\Pi:\, s\leq y/2} u_t\right)_{t\geq 0}$$

has the same distribution as the unconditioned process $(X_t)_{t\geq 0}$ with $X_0 = y$ (see [15]). The factor of 2 differs from [15] and arises from our choice of time scale. Therefore,

$$\mathbb{P}\{X_t > 0 | X_0 = y\} = \mathbb{P}\{\exists (s,u) \in \Pi : s \leq y/2 \text{ and } \gamma(u) > t\}$$
$$= 1 - \exp(-y\mu((t,\infty))/2),$$

and comparing with (5.1) gives $\mu((t,\infty)) = 2/t$.

Now take $\beta \in (0,1)$. It is shown in [16] that the mass thrown off the immortal lineage is determined by the jumps of a stable subordinator: if $M_s$ is the amount of mass thrown off during the time interval $[0,s]$, then

$$\mathbb{E}[e^{-\theta M_s}] = \exp\left(-s\frac{1+\beta}{\beta}\theta^\beta\right)$$
$$= \exp\left(-s\int_0^\infty (1 - e^{-\theta x})\nu(dx)\right),$$

where $\nu(dx) = \frac{1+\beta}{\Gamma(1-\beta)} x^{-(1+\beta)}\, dx$ is the Lévy measure of the subordinator.

Since the jump size of the subordinator corresponds to the initial size of the new family, the lifetime measure $\mu$ is given by

$$\mu((t,\infty)) = \int_0^\infty \mathbb{P}\{X_t > 0 | X_0 = x\}\nu(dx),$$

and so, from the above and an integration by parts,

$$\mu((t,\infty)) = \int_0^\infty (1 - e^{-x/t^{1/\beta}})\nu(dx)$$
$$= \frac{1+\beta}{\beta t}. \qquad \square$$

PROOF OF COROLLARY 1.1. Parts (a) to (d) follow immediately from Theorem 1.1. Part (e) will also follow from parts (e) and (f) of Theorem 1.1 after the following time and space change.

Define a new time parameter $u$ by $t = e^u$. If the MRCA at time $t$ lived at time $t - x$ on the original scale, then on the new time scale she lived at time $u - y$ where $t - x = e^{u-y}$. Solving for $y$, the MRCA age process in the new time scale is the process $(B_t)_{t \geq 0}$ given by $B_u = -\log(1 - e^{-u}A_{e^u})$. The process $(B_t)_{t \geq 0}$ is obtained by applying the construction (1.1) to the point process given by

$$\left\{\left(\log s, \log\left(1 + \frac{x}{s}\right)\right), (s,x) \in \Lambda\right\},$$

which is a Poisson point process on $\mathbb{R} \times \mathbb{R}_+$ with intensity measure $\lambda \otimes \rho$ where

$$\rho((y,\infty)) = \frac{1+\beta}{\beta(e^y - 1)}, \qquad y \in \mathbb{R}_{++}.$$

Note that, in general, a time and space change of the Poisson process $\Lambda$ gives a new Poisson point process, but the resulting intensity measure will not typically be of the form $\lambda \otimes \kappa$ for some measure $\kappa$: it is a special feature of $\mu$ and the transformation that the product measure structure is maintained in this case.

It is straightforward to check parts (e) and (f) of Theorem 1.1 that $(B_t)_{t \geq 0}$ has the stationary distribution

$$\pi(dx) = \frac{1+\beta}{\beta}e^{-x}(1 - e^{-x})^{1/\beta}\,dx$$

and that the distribution of $B_t$ converges to $\pi$ in total variation as $t \to \infty$. Part (e) of the corollary then follows from the observation that $\frac{A_t}{t} = 1 - e^{-B_{\log(t)}}$ and an elementary change of variables. $\square$

**6. An analogue of the Bessel-squared family.** In this last section, we determine the Laplace transform of the $(1+\beta)$-stable continuous state branching process conditioned on nonextinction. For $\beta = 1$, the unconditioned process is the Bessel-squared process with dimension 0 and it is well known that the conditioned process is the Bessel-squared process with dimension

4. A comparison of Laplace transforms will suggest that the unconditioned and conditioned processes for a given general $\beta$ may also be embedded in a family of Markov processes that is analogous to the Bessel-squared family. We show that such a family exists for each $\beta$ and we establish some of the properties of these families.

Recall that the transition probabilities of the unconditioned $(1 + \beta)$-stable continuous state branching process, $(X_t)_{t\geq 0}$, are characterized by the Laplace transforms,

$$\mathbb{E}^x[\exp(-\theta X_t)] = \exp(-x\theta(t\theta^\beta + 1)^{-1/\beta}).$$

Hence the transition probabilities of the conditioned process $Y$ are characterized by the Laplace transforms,

$$\mathbb{E}^y[\exp(-\theta Y_t)] = \frac{1}{y}\mathbb{E}^y[\exp(-\theta X_t)X_t]$$

$$= -\frac{1}{y}\frac{\partial}{\partial \theta}\mathbb{E}^y[\exp(-\theta X_t)]$$

$$= \exp(-y\theta(t\theta^\beta + 1)^{-1/\beta})(t\theta^\beta + 1)^{-(\beta+1)/\beta},$$

thus establishing (1.2).

Recall also that if $\beta = 1$, then $(X_t)_{t\geq 0}$ and $(Y_t)_{t\geq 0}$ are, up to a constant multiple, the Bessel-squared processes with dimensions 0 and 4, respectively. The Bessel-squared process, $(Z_t)_{t\geq 0}$, with dimension $d$ (not necessarily integral) is (up to constants) the Markov process characterized by the Laplace transforms,

$$\mathbb{E}^z[\exp(-\theta Z_t)] = \exp(-z\theta(t\theta + 1)^{-1})(t\theta + 1)^{-d/2}.$$

This suggests that for $0 < \beta < 1$ and $\delta \geq 0$, there might be a semigroup, $(P_t)_{t\geq 0}$, such that

(6.1) $$P_t \exp(-\theta \cdot)(x) = \exp(-x\theta(t\theta^\beta + 1)^{-1/\beta})(t\theta^\beta + 1)^{-\delta}.$$

We first verify that, for a fixed value of $x$, the right-hand side of (6.1) is the Laplace transform of a probability distribution (as a function of $\theta$). We already know that

$$\exp(-x\theta(t\theta^\beta + 1)^{-1/\beta})$$

is the Laplace transform of a probability measure, so it suffices to show that

$$(t\theta^\beta + 1)^{-\delta}$$

is also a Laplace transform of a probability distribution. Let $(S_t)_{t\geq 0}$ be the $\beta$-stable subordinator starting from $S_0 = 0$ normalized so that

$$\mathbb{E}[\exp(-\theta S_t)] = \exp(-\theta^\beta t)$$

424 S. N. EVANS AND P. L. RALPHand let $(T_t)_{t \geq 0}$ be the gamma subordinator starting from $T_0 = 0$ normalized so that for $t > 0$,
$$\mathbb{P}\{T_t \in dy\} = \frac{y^{t-1}}{\Gamma(t)} \exp(-y)\, dy,$$
and hence
$$\mathbb{E}[\exp(-\theta T_t)] = (\theta + 1)^{-t}.$$
Then if $S$ and $T$ are independent,
$$\mathbb{E}[\exp(-\theta S_{tT_\delta})] = \mathbb{E}[\exp(-\theta^\beta t T_\delta)] = (t\theta^\beta + 1)^{-\delta}.$$

We next verify that $(P_t)_{t \geq 0}$ is a semigroup. Observe that
$$P_s P_t \exp(-\theta \cdot)(x)$$
$$= (t\theta^\beta + 1)^{-\delta} P_s \exp(-\theta(t\theta^\beta + 1)^{-1/\beta} \cdot)(x)$$
$$= (t\theta^\beta + 1)^{-\delta} \exp(-x\theta(t\theta^\beta + 1)^{-1/\beta}(s\theta^\beta(t\theta^\beta + 1)^{-1} + 1)^{-1/\beta})$$
$$\quad \times (s\theta^\beta(t\theta^\beta + 1)^{-1} + 1)^{-\delta}$$
$$= \exp(-x\theta((s+t)\theta^\beta + 1)^{-1/\beta})((s+t)\theta^\beta + 1)^{-\delta}$$
$$= P_{s+t} \exp(-\theta \cdot)(x).$$

It is clear that $\lim_{t \downarrow 0} P_t \exp(-\theta \cdot)(x) = \exp(-\theta x)$ and so $\lim_{t \downarrow 0} P_t f(x) = f(x)$ for $f \in C_0(\mathbb{R}_+)$. Also, $\lim_{y \to x} P_t \exp(-\theta \cdot)(y) = P_t \exp(-\theta \cdot)(x)$, and so $\lim_{y \to x} P_t f(y) = f(x)$ for $f \in C_0(\mathbb{R}_+)$. The standard Feller construction gives that there is a strong Markov process, $(Z_t)_{t \geq 0}$, with semigroup $(P_t)_{t \geq 0}$.

This family of Markov processes shares many features of the Bessel-squared family. For example, it follows for $a, b > 0$, that
$$\mathbb{E}^{az}[\exp(-\theta a^{-1} Z_{bt})] = \exp\left(-y\theta\left(bt\left(\frac{\theta}{a}\right)^\beta + 1\right)^{-1/\beta}\right)\left(bt\left(\frac{\theta}{a}\right)^\beta + 1\right)^{-\delta}.$$

Thus, the process $(b^{-1/\beta} Z_{bt})_{t \geq 0}$ is Markovian with the same transition probabilities as $Z$. Similarly, if $Z_0 = 0$, then the process $(e^{-t/\beta} Z_{e^t})_{t \in \mathbb{R}}$ is Markovian and stationary.

Furthermore, if $(Z'_t)_{t \geq 0}$ and $(Z''_t)_{t \geq 0}$ are two independent such processes with parameters $\delta'$ and $\delta''$, then the process $(Z'_t + Z''_t)_{t \geq 0}$ also belongs to the family and has parameter $\delta' + \delta''$.

In a forthcoming paper, we will present a more thorough study of this family along the lines of [27, 28].

**Acknowledgments.** We thank the referee and Associate Editor for very careful readings of the paper and a number of suggestions that contributed to its clarity.

# REFERENCES

[1] Boxma, O., Perry, D., Stadje, W. and Zacks, S. (2006). A Markovian growth-collapse model. *Adv. in Appl. Probab.* **38** 221–243. MR2213972

[2] Borovkov, K. and Vere-Jones, D. (2000). Explicit formulae for stationary distributions of stress release processes. *J. Appl. Probab.* **37** 315–321. MR1780992

[3] Costa, O. L. V. and Dufour, F. (2008). Stability and ergodicity of piecewise deterministic Markov processes. *SIAM J. Control Optim.* **47** 1053–1077. MR2385873

[4] Colombo, G. and Dai Pra, P. (2001). A class of piecewise deterministic Markov processes. *Markov Process. Related Fields* **7** 251–287. MR1856497

[5] Davis, M. H. A. (1984). Piecewise-deterministic Markov processes: A general class of nondiffusion stochastic models. *J. Roy. Statist. Soc. Ser. B* **46** 353–388. MR790622

[6] Davis, M. H. A. (1993). *Markov Models and Optimization. Monographs on Statistics and Applied Probability* **49**. Chapman and Hall, London. MR1283589

[7] Dufour, F. and Costa, O. L. V. (1999). Stability of piecewise-deterministic Markov processes. *SIAM J. Control Optim.* **37** 1483–1502. MR1710229

[8] Dumas, V., Guillemin, F. and Robert, P. (2002). A Markovian analysis of additive-increase multiplicative-decrease algorithms. *Adv. in Appl. Probab.* **34** 85–111. MR1895332

[9] Daley, D. J. and Vere-Jones, D. (1988). *An Introduction to the Theory of Point Processes*. Springer, New York. MR950166

[10] Eliazar, I. (2007). The $M/G/\infty$ system revisited: Finiteness, summability, long range dependence, and reverse engineering. *Queueing Syst.* **55** 71–82. MR2293566

[11] Evans, S. N. and Perkins, E. (1990). Measure-valued Markov branching processes conditioned on nonextinction. *Israel J. Math.* **71** 329–337. MR1088825

[12] Evans, S. N. and Pitman, J. (1998). Stationary Markov processes related to stable Ornstein–Uhlenbeck processes and the additive coalescent. *Stochastic Process. Appl.* **77** 175–185. MR1649003

[13] Evans, S. N., Pitman, J. and Winter, A. (2006). Rayleigh processes, real trees, and root growth with re-grafting. *Probab. Theory Related Fields* **134** 81–126. MR2221786

[14] Evans, S. N. (1992). The entrance space of a measure-valued Markov branching process conditioned on nonextinction. *Canad. Math. Bull.* **35** 70–74. MR1157466

[15] Evans, S. N. (1993). Two representations of a conditioned superprocess. *Proc. Roy. Soc. Edinburgh Sect. A* **123** 959–971. MR1249698

[16] Etheridge, A. M. and Williams, D. R. E. (2003). A decomposition of the $(1+\beta)$-superprocess conditioned on survival. *Proc. Roy. Soc. Edinburgh Sect. A* **133** 829–847. MR2006204

[17] Feller, W. (1951). Diffusion processes in genetics. In *Proc. Second Berkeley Symp. Math. Statist. Probab.* 227–246. Univ. California Press, Berkeley. MR0046022

[18] Fitzsimmons, P. J., Fristedt, B. and Shepp, L. A. (1985). The set of real numbers left uncovered by random covering intervals. *Z. Wahrsch. Verw. Gebiete* **70** 175–189. MR799145

[19] Ingman, M., Kaessmann, H., Paabo, S. and Gyllensten, U. (2000). Mitochondrial genome variation and the origin of modern humans. *Nature* **408** 708–713.

[20] Jacod, J. and Skorokhod, A. V. (1996). Jumping Markov processes. *Ann. Inst. H. Poincaré Probab. Statist.* **32** 11–67. MR1373726

[21] Kallenberg, O. (2000). An extension of the basic Palm measure correspondence. *Probab. Theory Related Fields* **117** 113–131. MR1759510

[22] KNIGHT, F. B. (1981). *Essentials of Brownian Motion and Diffusion. Mathematical Surveys* **18**. Amer. Math. Soc., Providence, RI. MR613983

[23] MANDELBROT, B. B. (1972). Renewal sets and random cutouts. *Z. Wahrsch. Verw. Gebiete* **22** 145–157. MR0309162

[24] MORGAN, F. (2009). *Geometric Measure Theory: A Beginner's Guide*, 4th ed. Elsevier/Academic Press, Amsterdam. MR2455580

[25] MEYN, S. P. and TWEEDIE, R. L. (1993). *Markov Chains and Stochastic Stability*. Springer, London. MR1287609

[26] PFAFFELHUBER, P. and WAKOLBINGER, A. (2006). The process of most recent common ancestors in an evolving coalescent. *Stochastic Process. Appl.* **116** 1836–1859. MR2307061

[27] PITMAN, J. and YOR, M. (1981). Bessel processes and infinitely divisible laws. In *Stochastic Integrals (Proc. Sympos., Univ. Durham, Durham, 1980). Lecture Notes in Math.* **851** 285–370. Springer, Berlin. MR620995

[28] PITMAN, J. and YOR, M. (1982). A decomposition of Bessel bridges. *Z. Wahrsch. Verw. Gebiete* **59** 425–457. MR656509

[29] SIMON, D. and DERRIDA, B. (2006). Evolution of the most recent common ancestor of a population with no selection. *J. Stat. Mech. Theory Exp.* P05002.

[30] SHEPP, L. A. (1972). Covering the line with random intervals. *Z. Wahrsch. Verw. Gebiete* **23** 163–170. MR0322923

[31] TAKÁCS, L. (1962). *Introduction to the Theory of Queues*. Oxford Univ. Press, New York. MR0133880

[32] TRAVERS, S. A. A., CLEWLEY, J. P., GLYNN, J. R., FINE, P. E. M., CRAMPIN, A. C., SIBANDE, F., MULAWA, D., MCINERNEY, J. O. and MCCORMACK, G. P. (2004). Timing and reconstruction of the most recent common ancestor of the subtype C clade of Human Immunodeficiency Virus type 1. *J. Virol.* **78** 10501–10506.

[33] THOMSON, R., PRITCHARD, J. K., SHEN, P., OEFNER, P. J. and FELDMAN, M. W. (2000). Recent common ancestry of human Y chromosomes: Evidence from DNA sequence data. *Proc. Nat. Acad. Sci. India Sect. A* **97** 7360–7365.

[34] VOLKMAN, S. K., BARRY, A. E., LYONS, E. J., NIELSEN, K. M., THOMAS, S. M., MEHEE, C., THAKORE, S. S., DAY, K. P., WIRTH, D. F. and HARTL, D. L. (2001). Recent origin of *Plasmodium falciparum* from a single progenitor. *Science* **293** 482–484.

DEPARTMENT OF STATISTICS #3860
UNIVERSITY OF CALIFORNIA AT BERKELEY
367 EVANS HALL
BERKELEY, CALIFORNIA 94720-3860
USA
E-MAIL: evans@stat.berkeley.edu
URL: http://www.stat.berkeley.edu/users/evans/

DEPARTMENT OF EVOLUTION AND ECOLOGY
UNIVERSITY OF CALIFORNIA AT DAVIS
ONE SHIELDS AVENUE
DAVIS, CALIFORNIA 95616
USA
E-MAIL: plralph@ucdavis.edu
URL: http://www.eve.ucdavis.edu/plralph/